\documentclass[12pt]{article}
\usepackage{amsmath,amsfonts,amssymb}
\usepackage{theorem}
\usepackage[usenames]{color}

\nonstopmode

\theorembodyfont{\rmfamily}
\newtheorem{theorem}{Theorem}[section]
\newtheorem{proposition}[theorem]{Proposition}
\newtheorem{lemma}[theorem]{Lemma}
\newtheorem{definition}[theorem]{Definition}

\newtheorem{example}[theorem]{Example}
\newtheorem{remark}[theorem]{Remark}

\newcommand{\2}{\color[cmyk]{1,1,1,0} 2 \color{black}}

\begin{document}
$\,$%\vspace{10mm}

\begin{center}
\textrm{\LARGE Finding Rigged Configurations From Paths}
\vspace{5mm}
\end{center}
\begin{flushright}
{\large Reiho Sakamoto}\\
\end{flushright}

\section{Introduction}
In this lecture note, we review reformulation of the
bijection of Kerov--Kirillov--Reshetikhin \cite{KKR}
(and extension due to \cite{KR},
see also \cite{KSS})
$$\phi :\mathrm{Paths}
\longmapsto\mathrm{Rigged\, Configurations}$$
in terms of the crystal bases theory \cite{Kas}
and its application to the periodic box-ball systems
following \cite{S2} and \cite{KS3}.

The bijection $\phi$ was originally introduced in
order to show the so-called $X=M$ formula
(see \cite{O,Sch} for reviews) by using
its statistic preserving property.
Recently, another application of the
bijection $\phi$ to the box-ball systems \cite{TS,T}
was found \cite{KOSTY}.
In this context, the bijection $\phi$ plays the role
of the inverse scattering formalism \cite{GGKM,AC}
for the box-ball systems.

The original definition of the bijection $\phi$
is described by purely combinatorial language
such as box adding or removing procedures.
Purpose of this note is to clarify what the
representation theoretic origin of the bijection $\phi$ is.
Motivated by the connection with the box-ball systems,
consider the following isomorphism under the
affine combinatorial $R$:
\begin{equation}
u_l[0]\otimes b_1[0]\otimes\cdots\otimes b_L[0]\simeq
b_1'[-E_{l,1}]\otimes\cdots\otimes
b_L'[-E_{l,L}]\otimes u_l'[E_l],
\end{equation}
where $u_l\in B_l$ is the highest element and
$u_l'\in B_l$,
$b_k,b_k'\in B_{\lambda_k}$.
This is nothing but time evolution of the box-ball systems
\cite{HHIKTT,FOY}.
Then, $E_l$ here is related to shape of the rigged configuration
(see Eq.(\ref{eq:E=Q})) and it is the conserved quantity
of the box-ball systems introduced in \cite{FOY}.
By using this property, we introduce a table
containing purely algebraic data
(local energy distribution) from which we can read off
which letter 2 of path belongs to which row of
the rigged configuration.
As the result, we can reconstruct the map $\phi$
by using purely representation theoretic procedure.
Recently, the formalism is extended to
general elements of tensor products of the
Kirillov--Reshetikhin crystals of
$\mathfrak{sl}_n$ type \cite{S3}.
Again, this is achieved by extension of Eq.(\ref{eq:E=Q}).
This shows that the formalism here is quite natural.

The plan of this note is as follows.
In Section \ref{sec:crystal}, we prepare
minimal foundations of crystal theory.
In Section \ref{sec:LED}, we define
the local energy distribution.
In Section \ref{sec:main}, we present
our main result.
In Section \ref{sec:pbbs},
we explain some of applications of our
formalism for the box-ball system
with periodic boundary condition \cite{YT,YYT}
along with review of the inverse scattering
formalism for the periodic box-ball system
\cite{KTT}.
In Appendix, we consider tensor product of
highest paths in terms of the rigged configurations.

\section{Combinatorial R and energy function}
\label{sec:crystal}
\paragraph{Crystal.}
Let $B_k$ be the crystal of $k$-fold symmetric powers of the
vector (or natural) representation of $U_q(\mathfrak{sl}_2)$.
As the set, it is
\begin{equation}
B_k=\{(x_1,x_2)\in\mathbb{Z}^2_{\geq 0}\,\vert\,
x_1+x_2=k\}.
\end{equation}
We usually identify elements of $B_k$ as the 
semi-standard Young tableaux
\begin{equation}
(x_1,x_2)=
\fbox{$\overbrace{1\cdots 1}^{x_1}\overbrace{2\cdots 2}^{x_2}$}\, ,
\end{equation}
i.e., the number of letters $i$ contained in a tableau is $x_i$.
For example, the highest element $u_l\in B_l$ is
$u_l=(l,0)=\fbox{$111\cdots 1$}\,$.

For two crystals $B_k$ and $B_l$ of $U_q(\mathfrak{sl}_2)$,
one can define the tensor product
$B_k\otimes B_l=\{b\otimes b'\mid b\in B_k,b'\in B_l\}$.
Then we have a unique isomorphism $R:B_k\otimes B_l
\stackrel{\sim}{\rightarrow}B_l\otimes B_k$, i.e. a unique map
which commutes with actions of the Kashiwara operators
$\tilde{e}_i$, $\tilde{f}_i$.
We call this map combinatorial $R$
and usually write the map $R$ simply by $\simeq$.
We call elements of tensor product of crystals {\it paths}.

\paragraph{Affinization.}
Consider the affinization of the crystal $B$ \cite{KMN}.
As the set, it is
\begin{equation}
\mathrm{Aff}(B)=\{b[d]\, |\, b\in B,\, d\in\mathbb{Z}\}.
\end{equation}
Integers $d$ of $b[d]$ are called modes.
For the tensor product
$b_1[d_1]\otimes b_2[d_2]\in
\mathrm{Aff}(B_{k})\otimes\mathrm{Aff}(B_l)$,
we can lift the above definition of the
combinatorial $R$ as follows:
\begin{equation}
b_1[d_1]\otimes b_2[d_2]\stackrel{R}{\simeq}
b_2'[d_2-H(b_1\otimes b_2)]\otimes
b_1'[d_1+H(b_1\otimes b_2)],
\end{equation}
where $b_1\otimes b_2\simeq b_2'\otimes b_1'$
is the combinatorial $R$ defined in the above.

\paragraph{Explicit expressions.}
There is
piecewise linear formula to obtain the combinatorial $R$
and the energy function \cite{HHIKTT}.
This is suitable for computer programming.
For the affine combinatorial
$R:x[d]\otimes y[e]\simeq
\tilde{y}[e-H(x\otimes y)]\otimes
\tilde{x}[d+H(x\otimes y)]$,
we have
\begin{eqnarray}
&&\tilde{x}_i=x_i+Q_i(x,y)-Q_{i-1}(x,y),\qquad
\tilde{y}_i=y_i+Q_{i-1}(x,y)-Q_i(x,y),\nonumber\\
&&H(x\otimes y)=Q_0(x,y),\nonumber\\
&&Q_i(x,y)=\min (x_{i+1},y_i),
\end{eqnarray}
where we have expressed
$x=(x_1,x_2)$, $y=(y_1,y_2)$,
$\tilde{x}=(\tilde{x}_1,\tilde{x}_2)$ and
$\tilde{y}=(\tilde{y}_1,\tilde{y}_2)$.
All indices $i$ should be considered as
$i\in\mathbb{Z}/2\mathbb{Z}$.
There is another graphical method due to Nakayashiki--Yamada
\cite{NY} (see \cite{Shimo} for generalizations).
It is useful when we are going to prove mathematical statements.

\section{Local energy distribution}\label{sec:LED}
We express the isomorphism
$a\otimes b_1\simeq b_1'\otimes a'$
(with the energy function
$e_1:=H(a\otimes b_1)$)
by the following vertex diagram:
\begin{center}
\unitlength 13pt
\begin{picture}(4,4)
\put(0,2.0){\line(1,0){3.2}}
\put(1.6,1.0){\line(0,1){2}}
\put(-0.6,1.8){$a$}
\put(1.4,0){$b_1'$}
\put(1.4,3.2){$b_1$}
\put(3.4,1.8){$a'$}
\put(0.7,2.2){$e_1$}
\put(4.1,1.7){.}
\end{picture}
\end{center}
\begin{definition}
(1) For a given path $b=b_1\otimes b_2\otimes\cdots\otimes b_L$,
we define local energy $E_{l,j}$ by
$E_{l,j}:=H(u_l^{(j-1)}\otimes b_j)$.
Here $u_l^{(j-1)}$ are defined in the
following diagram
(we set $u_l^{(0)}:=u_l$):
\begin{center}
\unitlength 13pt
\begin{picture}(22,5)(0,-0.5)
\multiput(0,0)(5.8,0){2}{
\put(0,2.0){\line(1,0){4}}
\put(2,0){\line(0,1){4}}
}
\put(-0.9,1.8){$u_l$}
\put(0.5,2.3){$E_{l,1}$}
\put(1.7,4.2){$b_1$}
\put(1.7,-0.8){$b_1'$}
\put(4.2,1.8){$u_l^{(1)}$}
\put(6.3,2.3){$E_{l,2}$}
\put(7.5,4.2){$b_2$}
\put(7.5,-0.8){$b_2'$}
\put(10.0,1.8){$u_l^{(2)}$}
\multiput(11.5,1.8)(0.3,0){10}{$\cdot$}
\put(14.7,1.8){$u_l^{(L-1)}$}
\put(17,0){
\put(0,2.0){\line(1,0){4}}
\put(2,0){\line(0,1){4}}
}
\put(17.4,2.3){$E_{l,L}$}
\put(18.7,4.2){$b_L$}
\put(18.7,-0.8){$b_L'$}
\put(21.2,1.8){$u_l^{(L)}$}
\end{picture}
\end{center}
We define $E_{0,j}=0$ for all $1\leq j\leq L$.
We also use the notation
$E_l:=\sum_{j=1}^L E_{l,j}$
which coincides with the conserved quantity in \cite{FOY}.\\
(2) We define operator $T_l$ by
$T_l(b)=b_1'\otimes b_2'\otimes\cdots\otimes b_L'$
(see the above diagram).
\hfill$\square$
\end{definition}

\begin{lemma}\label{s:lem:E=0,1}
For a given path $b=b_1\otimes b_2\otimes\cdots\otimes b_L$,
we have $E_{l,j}-E_{l-1,j}=0$ or 1,
for all $l>0$ and for all $1\leq j\leq L$.
\hfill$\square$
\end{lemma}
\begin{definition}
The local energy distribution
is a table containing $(E_{l,j}-E_{l-1,j}=0,1)$
at the position $(l,j)$, i.e., at the $l$ th row and
the $j$ th column.
\hfill$\square$
\end{definition}

\section{Results}\label{sec:main}
\subsection{Statement}
\begin{theorem}\label{th:main}
Let $b=b_1\otimes b_2\otimes\cdots\otimes b_L
\in B_{\lambda_1}\otimes B_{\lambda_2}\otimes\cdots\otimes
B_{\lambda_L}$ be an arbitrary path.
$b$ can be highest weight or non-highest weight.
Set $N=E_1(b)$.
We determine the pair of numbers
$(\mu_1,r_1)$, $(\mu_2,r_2)$, $\cdots$, $(\mu_N,r_N)$
by the following procedure from Step 1 to Step 4.
Then the resulting $(\lambda ,(\mu ,r))$ coincides with
the (unrestricted) rigged configuration $\phi (b)$.
\begin{enumerate}
\item
Draw the local energy distribution for $b$.
\item
Starting from the rightmost 1 in the $l=1$ st row,
pick one 1 from each successive row.
The one in the $(l+1)$ th row must be weakly right of
the one selected in the $l$ th row.
If there is no such 1 in the $(l+1)$ th row,
the position of the lastly picked 1 is called $(\mu_1,j_1)$.
Change all selected 1 into 0.
\item
Repeat Step 2 for $(N-1)$ times to further determine
$(\mu_2,j_2)$, $\cdots$, $(\mu_N,j_N)$
thereby making all 1 into 0.
\item
Determine $r_1,\cdots,r_N$ by
\begin{equation}\label{s:eq:rigging}
r_k=\sum_{i=1}^{j_k-1}\min (\mu_k,\lambda_i)
+E_{\mu_k,j_k}-2\sum_{i=1}^{j_k}E_{\mu_k,i}.
\end{equation}
\end{enumerate}
\end{theorem}
{\bf Sketch of proof.}
The key formula is
\begin{equation}\label{eq:E=Q}
E_l=\sum_{i=1}^N\min (\mu_i,l).
\end{equation}
The rest is combinatorial arguments whose details we
left to \cite{S2}.
\hfill$\square$

In order to depict the rigged configurations,
we usually use Young diagrammatic expression whose
rows have lengths $\mu_1,\cdots,\mu_N$,
and we put integers $r_1,\cdots,r_N$ on the right
of rows $\mu_1,\cdots,\mu_N$, respectively.
Integers $r_k$ are called {\it riggings} corresponding to $\mu_k$.
The letters 1 in the local energy distribution precisely records
combinatorial procedure of $\phi$, i.e.,
letter 1 at $l$ th row, $k$ th column of the local energy
distribution corresponds to box addition to $l$ th
column of some row of rigged configuration.

The groups $\mu_i$ obtained here represent solitons
(moving at velocity $\mu_i$ with respect to $T_\infty$)
contained in a path.
In fact, we have the following property \cite{KOSTY}.
\begin{proposition}
Given the (unrestricted) rigged configuration corresponding
to $b$:
\begin{equation}\label{s:prop:T_l:1}
b\xrightarrow{\,\phi\,}
\bigl( \lambda
,(\mu_j ,r_j)_{j=1}^N\bigl).
\end{equation}
Then, corresponding to $T_l(b)$, we have
\begin{equation}\label{s:prop:T_l:2}
T_l(b)\xrightarrow{\,\phi\,}
\bigl( \lambda
,(\mu_j ,r_j+\min (\mu_j,l))_{j=1}^N\bigl).
\end{equation}
\hfill$\square$
\end{proposition}%\bigskip
For the proof,
see Proposition 2.6 of \cite{KOSTY}
(see also Proposition 2.3 of \cite{S2}).

\begin{remark}
A crystal interpretation of the inverse map $\phi^{-1}$ is known
\cite{S1} for the cases
$B^{1,s_1}\otimes\cdots\otimes B^{1,s_L}$.
This formalism gives recursive description of $\phi^{-1}$
with respect to rank of $\mathfrak{sl}_n$ and,
owing to this property,
it is substantially used in \cite{KSY}.
However, our formalism seems to have different origin
from \cite{S1},
since it can be generalized \cite{S3}
to wider class of representations
$B^{r_1,s_1}\otimes\cdots\otimes B^{r_L,s_L}$.
In this generalization, the procedure is done almost
separately with respect to the rank.
\hfill$\square$
\end{remark}

\subsection{Example}
Consider the following path:
\begin{equation}
b=
\fbox{1111}\otimes\fbox{11}\otimes\fbox{2}\otimes
\fbox{1122}\otimes\fbox{1222}\otimes\fbox{1}\otimes
\fbox{2}\otimes\fbox{122}
\end{equation}
Corresponding to Step 1,
the local energy distribution is given by
the following table
($j$ stands for column coordinate of the table).
\begin{center}
\begin{tabular}{l|cccccccc}
                 &1111&11&2  &1122 &1222 &1&2 &122  \\\hline
$E_{1,j}-E_{0,j}$&0   &0 &1  &0    &1    &0&1 &0    \\
$E_{2,j}-E_{1,j}$&0   &0 &0  &1    &0    &0&0 &1    \\
$E_{3,j}-E_{2,j}$&0   &0 &0  &1    &0    &0&0 &0    \\
$E_{4,j}-E_{3,j}$&0   &0 &0  &0    &1    &0&0 &0    \\
$E_{5,j}-E_{4,j}$&0   &0 &0  &0    &1    &0&0 &0    \\
$E_{6,j}-E_{5,j}$&0   &0 &0  &0    &0    &0&0 &1    \\
$E_{7,j}-E_{6,j}$&0   &0 &0  &0    &0    &0&0 &0    \\
\end{tabular}
\end{center}
Following Step 2 and Step 3, letters 1
contained in the above table are found
to be classified into 3 groups,
as indicated in the following table.
\begin{center}
\begin{tabular}{l|cccccccc}
                 &1111&11&2 &1122&1222      &1  &2 &122       \\\hline
$E_{1,j}-E_{0,j}$&    &  &3 &    &$2^{\ast}$&   &1 &          \\
$E_{2,j}-E_{1,j}$&    &  &  &3   &          &   &  &$1^{\ast}$\\
$E_{3,j}-E_{2,j}$&    &  &  &3   &          &   &  &          \\
$E_{4,j}-E_{3,j}$&    &  &  &    &3         &   &  &          \\
$E_{5,j}-E_{4,j}$&    &  &  &    &3         &   &  &          \\
$E_{6,j}-E_{5,j}$&    &  &  &    &          &   &  &$3^{\ast}$\\
$E_{7,j}-E_{6,j}$&    &  &  &    &          &   &  &          \\
\end{tabular}
\end{center}
{}From the above table, we see that
the cardinalities of groups 1, 2 and 3 are
2, 1 and 6, respectively.
Also, in the above table,
positions of $(\mu_1,j_1)$, $(\mu_2,j_2)$
and $(\mu_3,j_3)$ are indicated by asterisks.
Their explicit locations are $(\mu_1,j_1)=(2,8)$,
$(\mu_2,j_2)=(1,5)$
and $(\mu_3,j_3)=(6,8)$.

Now we evaluate riggings $r_i$ according to
Eq.(\ref{s:eq:rigging}).
\begin{eqnarray*}
r_1&=&\sum_{i=1}^{8-1} \min(2,\lambda_i)
      +E_{2,8}-2\sum_{i=1}^8E_{2,i}\\
   &=&(2+2+1+2+2+1+1)+1-2(0+0+1+1+1+0+1+1)\\
   &=&2,\\
r_2&=&\sum_{i=1}^{5-1} \min(1,\lambda_i)
      +E_{1,5}-2\sum_{i=1}^5E_{1,i}\\
   &=&(1+1+1+1)+1-2(0+0+1+0+1)\\
   &=&1,\\
%\end{eqnarray*}
%\begin{eqnarray*}
r_3&=&\sum_{i=1}^{8-1} \min(6,\lambda_i)
      +E_{6,8}-2\sum_{i=1}^8E_{6,i}\\
   &=&(4+2+1+4+4+1+1)+2-2(0+0+1+2+3+0+1+2)\\
   &=&1.
\end{eqnarray*}
Therefore we obtain
$(\mu_1,r_1)=(2,2)$, $(\mu_2,r_2)=(1,1)$ and $(\mu_3,r_3)=(6,1)$.
This coincides with the calculation based on the original
combinatorial definition of the map $\phi$.

\section{Application to periodic box-ball system}
\label{sec:pbbs}
\subsection{Definition}

In this section, we consider application of Theorem \ref{th:main}
to the periodic box-ball system (pBBS).
Many part of this section is contained in \cite{KTT,KS1}.
We exclusively treat $\mathfrak{sl}_2$ type path $b$ of the form
$b\in B_1^{\otimes L}$.
The pBBS is the BBS with periodic boundary condition
and its definition rely on the following fact.
\begin{proposition}\label{prop:v_l}
Define $v_l\in B_l$ by
\begin{equation}\label{def:v_l}
u_l\otimes b\stackrel{R}{\simeq}
T_l(b)\otimes v_l.
\end{equation}
Then we have
\begin{equation}\label{def:barT}
v_l\otimes b\stackrel{R}{\simeq}
b^\ast\otimes v_l,
\end{equation}
where $b^\ast\in B_1^{\otimes L}$.
\hfill$\square$
\end{proposition}
For the proof,
see Proposition 2.1 of \cite{KTT} and the comment following it.

\begin{definition}
We define operator of pBBS $\bar{T}_l$ by
$\bar{T}_l(b)=b^\ast\in B_1^{\otimes L}$,
where $b^\ast$ is obtained in the right hand side
of Eq.(\ref{def:barT}).
\hfill$\square$
\end{definition}
Note that $\bar{T}_1$ is merely the cyclic shift operator on a path.

\begin{example}
The time evolutions $b, \bar{T}_l(b), \ldots, \bar{T}^9_l(b)$ 
of the state $b$ on the top line are listed downward
for $l=2$ and $3$. The system size is  $L=14$.
We omit the symbol $\otimes$
and frames of tableaux.
\[
%\begin{array}{lllllllllllllllllllllllllllll}
\begin{array}{llll}
\quad \;\hbox{evolution under } \bar{T}_2 
\quad \qquad \quad\qquad \hbox{evolution under } \bar{T}_3 \\
1\; 1\; \2\; 1\; 1\; \2\; \2\; 1\; 1\; 1\; 1\; \2\; \2\; \2\;
\qquad\quad\; 1\; 1\; \2\; 1\; 1\; \2\; \2\; 1\; 1\; 1\; 1\; \2\; \2\; \2\\
\2\; \2\; 1\; \2\; 1\; 1\; 1\; \2\; \2\; 1\; 1\; 1\; 1\; \2\;
\qquad\quad\; \2\; \2\; 1\; \2\; \2\; 1\; 1\; \2\; \2\; 1\; 1\; 1\; 1\; 1\\
1\; \2\; \2\; 1\; \2\; \2\; 1\; 1\; 1\; \2\; \2\; 1\; 1\; 1\;
\qquad\quad\; 1\; 1\; \2\; 1\; 1\; \2\; \2\; 1\; 1\; \2\; \2\; \2\; 1\; 1\\
1\; 1\; 1\; \2\; 1\; \2\; \2\; \2\; 1\; 1\; 1\; \2\; \2\; 1\;
\qquad\quad\; \2\; 1\; 1\; \2\; 1\; 1\; 1\; \2\; \2\; 1\; 1\; 1\; \2\; \2\\
\2\; 1\; 1\; 1\; \2\; 1\; 1\; \2\; \2\; \2\; 1\; 1\; 1\; \2\;
\qquad\quad\; 1\; \2\; \2\; 1\; \2\; \2\; 1\; 1\; 1\; \2\; \2\; 1\; 1\; 1\\
1\; \2\; \2\; 1\; 1\; \2\; 1\; 1\; 1\; \2\; \2\; \2\; 1\; 1\;
\qquad\quad\; 1\; 1\; 1\; \2\; 1\; 1\; \2\; \2\; \2\; 1\; 1\; \2\; \2\; 1\\
1\; 1\; 1\; \2\; \2\; 1\; \2\; 1\; 1\; 1\; 1\; \2\; \2\; \2\;
\qquad\quad\; \2\; \2\; 1\; 1\; \2\; 1\; 1\; 1\; 1\; \2\; \2\; 1\; 1\; \2\\
\2\; \2\; 1\; 1\; 1\; \2\; 1\; \2\; \2\; 1\; 1\; 1\; 1\; \2\;
\qquad\quad\; 1\; 1\; \2\; \2\; 1\; \2\; \2\; 1\; 1\; 1\; 1\; \2\; \2\; 1\\
1\; \2\; \2\; \2\; 1\; 1\; \2\; 1\; 1\; \2\; \2\; 1\; 1\; 1\;
\qquad\quad\; \2\; 1\; 1\; 1\; \2\; 1\; 1\; \2\; \2\; \2\; 1\; 1\; 1\; \2\\
1\; 1\; 1\; \2\; \2\; \2\; 1\; \2\; 1\; 1\; 1\; \2\; \2\; 1\;
\qquad\quad\; 1\; \2\; \2\; 1\; 1\; \2\; 1\; 1\; 1\; 1\; \2\; \2\; \2\; 1
\end{array}
\]
There are three solitons with amplitudes $3,2$ and $1$ 
traveling to the right.
\hfill$\square$
\end{example}

\subsection{Basic procedures}
In the rest of the note, we exclusively consider the path
$b\in B_1^{\otimes L}$ where number of \fbox{2} is equal to or less than
that of \fbox{1}.
The other case follows from this case by virtue of
Proposition 2.3 of \cite{KTT}.
In order to analyze the pBBS by using Theorem \ref{th:main},
we follow the following procedures.

\begin{enumerate}
\item
Instead of using $u_l$,
use $v_l$ of Eq.(\ref{def:v_l})
to calculate the local energy distribution.
We express energy function appearing here
as $\bar{E}_{l,j}$ and $\bar{E}_l=\sum_j\bar{E}_{l,j}$.
See the following diagram:
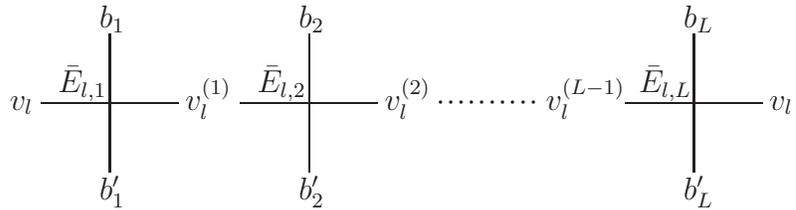
\begin{center}
\unitlength 13pt
\begin{picture}(22,5)(0,-0.5)
\multiput(0,0)(5.8,0){2}{
\put(0,2.0){\line(1,0){4}}
\put(2,0){\line(0,1){4}}
}
\put(-0.9,1.8){$v_l$}
\put(0.5,2.3){$\bar{E}_{l,1}$}
\put(1.7,4.2){$b_1$}
\put(1.7,-0.8){$b_1'$}
\put(4.2,1.8){$v_l^{(1)}$}
\put(6.3,2.3){$\bar{E}_{l,2}$}
\put(7.5,4.2){$b_2$}
\put(7.5,-0.8){$b_2'$}
\put(10.0,1.8){$v_l^{(2)}$}
\multiput(11.5,1.8)(0.3,0){10}{$\cdot$}
\put(14.7,1.8){$v_l^{(L-1)}$}
\put(17,0){
\put(0,2.0){\line(1,0){4}}
\put(2,0){\line(0,1){4}}
}
\put(17.4,2.3){$\bar{E}_{l,L}$}
\put(18.7,4.2){$b_L$}
\put(18.7,-0.8){$b_L'$}
\put(21.2,1.8){$v_l$}
\end{picture}
\end{center}

\item
Pick one of the lowest 1.
Pick 1 in $l$th row which is weakly
left of the already selected 1 in $(l+1)$th row.
If there is no such 1, return to the rightmost column
and search 1.

\item
Find a vertical line such that no soliton cross the line.
We call the line {\it seam}.
By using cyclic shift $\bar{T}_1$,
move the seam to the left end of the path.
Write this procedure as $b_+=\bar{T}_1^d(b)$.

\item
Apply Theorem \ref{th:main} to $b_+$.
(In order to obtain the local energy distribution here,
we only have to rotate periodic version of the local energy
distribution obtained in Step 1 according to $\bar{T}_1^d$).
\end{enumerate}

We can always find seam of a path due to the
following simple property
and Proposition \ref{prop:highest*highest}.
\begin{lemma}
For arbitrary element $b\in B_1^{\otimes L}$,
there exists integer $d$ such that
$b_+=\bar{T}_1^d(b)$ is highest weight.
\hfill $\square$
\end{lemma}
This assertion is proved by
elementary argument.
See Example 3.2 of \cite{KTT}.

\bigskip

The $b_+$ in this lemma can be used in Step 3 of the above procedure.
Consider the path $b_+^{\otimes n}$ and calculate the local energy
distribution using $u_l$ (not $v_l$).
Combining Eq.(\ref{def:barT}) and
the argument used in Proposition \ref{prop:highest*highest},
we can show that the local energy
distribution is $n$ times repetition of that of $b_+$.

\subsection{Action variables}
\begin{definition}
Given $b\in B_1^{\otimes L}$,
choose the highest element $b_+$ such that there exists integer $d$
such that $b=\bar{T}_1(b_+)$.
Apply the bijection $\phi$ and obtain
$\phi (b_+)=\left((1^L),(\mu,J)\right)$.
Then $\mu$ is called {\it action variable} of $b$.
\hfill$\square$
\end{definition}
\begin{proposition}
For any $l\in\mathbb{Z}_{\geq 1}$,
the action variable of $\bar{T}_l(b)$ is
equal to that of $b$.
\end{proposition}
{\bf Sketch of proof.}
Follows from the relation $\bar{E}_k(\bar{T}_l(b))=\bar{E}_k(b)$
which is the consequence of the Yang--Baxter relation
for the affine crystals.
See Theorem 2.2 of \cite{KTT} for more details.
\hfill$\square$

\subsection{Definition of the angle variables}
Thanks to Proposition \ref{prop:v_l},
our basic strategy is to
embed suitably cut periodic path $b$
into the usual infinite system as $b\otimes b\otimes
\cdots\otimes b$.
For this purpose, it is convenient to use
the highest path obtained by $b_+=\bar{T}_1^d(b)$
with suitable $d$.
However, this correspondence between $b$ and $(d,b_+)$
is not unique in general.
In the following, we give prescriptions to cope with
this ambiguity.

\subsubsection{Notations}
We fix some notations used in the following arguments.
Let $b_+$ be a highest element of $B_1^{\otimes L}$
and the corresponding rigged configuration be
\begin{equation}
b_+\xrightarrow{\,\phi\,}
\bigl(
(1^L),(\mu_i,J_i)_{i=1}^N
\bigl).
\end{equation}
Denote the multiplicity of $k$ in $(\mu_i)_{i=1}^N$
by $m_k$, and the riggings corresponding to
length $k$ rows by $J^{(k)}_1\leq J^{(k)}_2\leq
\cdots\leq J^{(k)}_{m_k}$.
Let the distinct lengths of rows of $(\mu_i)$ be
$k_1<k_2<\cdots <k_s$.
We denote the set of distinct lengths of rows of $(\mu_i)$
as $H=\{k_1,\ldots,k_s\}$.
Finally, we define the set of all possible riggings
as follows:
\begin{equation}
\mathrm{Rig}_L(\mu)=
\left\{
\left(
J_i^{(k)}
\right)_{1\leq i\leq m_k,k\in H}\in
\mathbb{Z}^{m_{k_1}}\times\cdots\times
\mathbb{Z}^{m_{k_s}}
\,\biggl|\,
0\leq J_1^{(k)}\leq\cdots\leq
J_{m_k}^{(k)}\leq p_k
\right\}.
\end{equation}
We sometimes omit $L$ of $\mathrm{Rig}_L(\mu)$ such as
$\mathrm{Rig}(\mu)$.
Here integer $p_k$ is called the vacancy number defined by
\begin{equation}
p_k=L-2\sum_{i=1}^N\min (k,\mu_i).
\end{equation}
In the present setting, we have $0\leq p_{k_s}<p_{k_{s-1}}
<\cdots <p_{k_1}$
(positivity $0\leq p_{k_i}$ follows from the fact that
$b_+$ is highest weight, and other inequalities $<$
follow from the shape of the quantum space $(1^L)$).

\subsubsection{Extension of riggings}
First we give motivations for extension of riggings.
{}From Proposition \ref{prop:highest*highest},
the rigged configuration corresponding to
$b_+^{\otimes n}$ have $n\times m_k$ rows of length $k$,
and the associated riggings takes the form
\begin{align}
&J^{(k)}_1\leq J^{(k)}_2\leq\cdots\leq J^{(k)}_{m_k}
\nonumber\\
\leq\,&
J^{(k)}_1+p_k\leq J^{(k)}_2+p_k\leq\cdots\leq J^{(k)}_{m_k}+p_k
\nonumber\\
\leq\,&\cdots
\nonumber\\
\leq\,&
J^{(k)}_1+(n-1)p_k\leq J^{(k)}_2+(n-1)p_k\leq\cdots\leq
J^{(k)}_{m_k}+(n-1)p_k.
\end{align}

In view of this observation, we define extension
of riggings
\begin{equation}
\iota :\left( J^{(k)}_i\right)_{1\leq i\leq m_k}
\longmapsto
\left( J^{(k)}_i\right)_{i\in\mathbb{Z}}
\end{equation}
by the relation
\begin{equation}
J^{(k)}_{i+m_k}=J^{(k)}_i+p_k.
\end{equation}
This $\iota (J)$ can be considered as an element
of the following set
\begin{equation}
\bar{\mathcal{J}}_k=
\left\{
\left( J^{(k)}_i\right)_{i\in\mathbb{Z}}
\Bigl|
J^{(k)}_i\in\mathbb{Z},\,
J^{(k)}_i\leq J^{(k)}_{i+1},\,
J^{(k)}_{i+m_k}=J^{(k)}_i+p_k,\,
\forall i
\right\} .
\end{equation}
We also define
\begin{equation}
\bar{\mathcal{J}}=\bar{\mathcal{J}}(\mu)=
\bar{\mathcal{J}}_{k_1}\times
\bar{\mathcal{J}}_{k_2}\times\cdots\times
\bar{\mathcal{J}}_{k_s}.
\end{equation}

\subsubsection{Slide $\sigma_l$ and equivalence relation}
On the extended riggings, we define the following
important operations.

\begin{definition}
For $l\in\mathbb{Z}_{\geq 1}$,
we define $\sigma_l:\bar{\mathcal{J}}_k
\mapsto\bar{\mathcal{J}}_k$
by
\begin{equation}\label{def:slide}
\sigma_l:
\left(
J^{(k)}_i
\right)_{i\in\mathbb{Z}}\longmapsto
\left(
J^{(k)}_{i+\delta_{l,k}}+2\min (l,k)
\right)_{i\in\mathbb{Z}}.
\end{equation}
We define abelian group $\mathcal{A}$
by
\begin{equation}
\mathcal{A}=
\left\{
\sigma_{k_1}^{n_1}\sigma_{k_2}^{n_2}\cdots\sigma_{k_s}^{n_s}
\bigl|
n_1,n_2,\ldots,n_s\in\mathbb{Z}
\right\}.
\end{equation}
We call an element of $\mathcal{A}$ {\it slide}.
\hfill$\square$
\end{definition}
We naturally define $\sigma_l$ on $\bar{\mathcal{J}}$
by $\sigma_l(\bar{\mathcal{J}})
=\sigma_l(\bar{\mathcal{J}}_{k_1})\times
\sigma_l(\bar{\mathcal{J}}_{k_2})\times\cdots\times
\sigma_l(\bar{\mathcal{J}}_{k_s})$.
\begin{definition}
We define equivalence relation $\simeq$ between
$J,K\in\bar{\mathcal{J}}$ by the following condition:
$J\simeq K$ if $\exists\sigma\in\mathcal{A}$
such that $J=\sigma(K)$.
\hfill$\square$
\end{definition}

We have the following standard form
with respect to the above $\simeq$.
\begin{proposition}\label{prop:standard-form}
For any $\bar{J}\in\bar{\mathcal{J}}(\mu)$,
there exist $d\in\mathbb{Z}$ and $J\in\mathrm{Rig}(\mu)$
such that $\bar{J}\simeq\iota (J)+d$.
\end{proposition}
{\bf Sketch of proof.}
There is a general algorithm to derive the standard form.
Basis of the algorithm is the relation $p_{k_i}>0$ ($2\leq i\leq s$).
See Lemma 3.9 of \cite{KTT} for more details.\bigskip
\hfill$\square$

\subsubsection{Interpretation of $\sigma_l$}
By using Theorem \ref{th:main}, we can give interpretation of
slide $\sigma_l$ in terms of the bijection $\phi$.
Suppose we have two expressions $p=\bar{T}_1^{d}(p_+)$ and
$p=\bar{T}_1^{d'}(p_+')$ with
highest paths $p_+$ and $p_+'$.
Let the rigged configuration corresponding to $p_+$
(resp. $p_+'$) be $(\mu,J)$ (resp. $(\mu,J')$).
Draw local energy distribution of $p$.
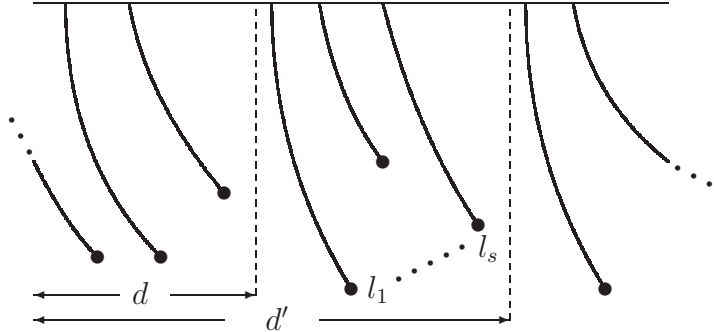
\begin{center}
\unitlength 12pt
\begin{picture}(20,10)
\put(0,10){\line(1,0){20}}
\multiput(7,0.8)(0,0.4){23}{\line(0,1){0.2}}
\multiput(15,0)(0,0.4){25}{\line(0,1){0.2}}
\put(2.7,0.8){\vector(-1,0){2.7}}
\put(4.3,0.8){\vector(1,0){2.7}}
\put(3.1,0.5){$d$}
\put(6,0){\vector(-1,0){6}}
\put(9,0){\vector(1,0){6}}
\put(7.3,-0.3){$d'$}
\thicklines
\qbezier(1,10)(1,5)(4,2)
\put(4,2){\circle*{0.4}}
\qbezier(3,10)(3.5,7)(6,4)
\put(6,4){\circle*{0.4}}
\qbezier(0,5)(1,3)(2,2)
\put(2,2){\circle*{0.4}}
\qbezier(7.5,10)(7.5,5)(10,1)
\put(10,1){\circle*{0.4}}
\put(10.5,0.7){$l_1$}
\qbezier(9,10)(9.5,7)(11,5)
\put(11,5){\circle*{0.4}}
\qbezier(11,10)(12,6)(14,3)
\put(14,3){\circle*{0.4}}
\put(14,2){$l_s$}
\multiput(11.5,1.3)(0.5,0.25){5}{\circle*{0.2}}
\qbezier(15.5,10)(15.5,5)(18,1)
\put(18,1){\circle*{0.4}}
\qbezier(17,10)(17.5,7)(20,5)
\multiput(20.3,4.8)(0.5,-0.25){3}{\circle*{0.2}}
\multiput(-0.7,6.3)(0.28,-0.5){3}{\circle*{0.2}}
\end{picture}
\end{center}
If the difference between $p_+$ and $p_+'$ is
solitons of lengths $l_1,\ldots,l_s$,
then we have
\begin{equation}\label{eq:rel_sigma}
\iota (J')+d'=\prod_{i=1}^s\sigma_{l_i}\iota (J)+d.
\end{equation}
Hence $\iota (J')+d'\simeq \iota (J)+d$.
This follows from direct calculation using Eq.(\ref{s:eq:rigging}).

\begin{example}
Consider the path
$p=221221112221111$
and $p'=111222111122122$ (omitting $\otimes$ 
and frames of tableaux, $p=\bar{T}_1^5(p')$).
Then the local energy distribution takes the following form
(all letters 0 are suppressed):
\begin{center}
\begin{tabular}{l|ccccccccccccccc}
                 &2 &2 &1 &2 &2 &1 &1 &1 &2 &2 &2 &1 &1 &1 &1 \\\hline
$\bar{E}_{1,j}-\bar{E}_{0,j}$&1 &  &  &1 &  &  &  &  &1 &  &  &  &  &  &  \\
$\bar{E}_{2,j}-\bar{E}_{1,j}$&  &1 &  &  &  &  &  &  &  &1 &  &  &  &  &  \\
$\bar{E}_{3,j}-\bar{E}_{2,j}$&  &  &  &  &1 &  &  &  &  &  &1 &  &  &  &  \\
\end{tabular}
\end{center}
We are setting $d=0$, $d'=5$, $l_1=3$ and $l_2=1$
(vacancy numbers are $p_3=1$ and $p_1=9$).
Here rigged configurations are
\begin{center}
\unitlength 13pt
\begin{picture}(15,4.5)
\put(0,0){\line(1,0){1}}
\multiput(0,0)(1,0){2}{\line(0,1){3}}
\multiput(0,1)(0,1){3}{\line(1,0){3}}
\multiput(2,1)(1,0){2}{\line(0,1){2}}
\put(3.2,2.2){$-3$}
\put(3.2,1.2){$-3$}
\put(1.2,0.2){$0$}
\put(1.3,3.8){$p$}
\put(10,0){
\put(0,0){\line(1,0){1}}
\multiput(0,0)(1,0){2}{\line(0,1){3}}
\multiput(0,1)(0,1){3}{\line(1,0){3}}
\multiput(2,1)(1,0){2}{\line(0,1){2}}
\put(3.2,2.2){$1$}
\put(3.2,1.2){$0$}
\put(1.2,0.2){$8$}
\put(1.3,3.8){$p'$}
}
\end{picture}
\end{center}
Calculation goes as follows:

\begin{center}
\unitlength 13pt
\begin{picture}(25,4)
\put(0,0){\line(1,0){1}}
\multiput(0,0)(1,0){2}{\line(0,1){3}}
\multiput(0,1)(0,1){3}{\line(1,0){3}}
\multiput(2,1)(1,0){2}{\line(0,1){2}}
\put(3.2,2.2){$-3$}
\put(3.2,1.2){$-3$}
\put(1.2,0.2){$0$}
\put(4,1.7){$\quad\xrightarrow{\sigma_1}$}
\put(7,0){
\put(0,0){\line(1,0){1}}
\multiput(0,0)(1,0){2}{\line(0,1){3}}
\multiput(0,1)(0,1){3}{\line(1,0){3}}
\multiput(2,1)(1,0){2}{\line(0,1){2}}
\put(3.2,2.2){$-1$}
\put(3.2,1.2){$-1$}
\put(1.2,0.2){$11$}
\put(4,1.7){$\quad\xrightarrow{\sigma_3}$}
}
\put(14,0){
\put(0,0){\line(1,0){1}}
\multiput(0,0)(1,0){2}{\line(0,1){3}}
\multiput(0,1)(0,1){3}{\line(1,0){3}}
\multiput(2,1)(1,0){2}{\line(0,1){2}}
\put(3.2,2.2){$6$}
\put(3.2,1.2){$5$}
\put(1.2,0.2){$13$}
\put(4.5,1.7){$=$}
}
\put(22,0){
\put(0,0){\line(1,0){1}}
\multiput(0,0)(1,0){2}{\line(0,1){3}}
\multiput(0,1)(0,1){3}{\line(1,0){3}}
\multiput(2,1)(1,0){2}{\line(0,1){2}}
\put(3.2,2.2){$1$}
\put(3.2,1.2){$0$}
\put(1.2,0.2){$8$}
\put(-2,1.7){$5\,\, +$}
}
\end{picture}
\end{center}
This coincides with
Eq.(\ref{eq:rel_sigma}).
\hfill$\square$
\end{example}

\subsubsection{Angle variables and inverse scattering formalism}
\begin{definition}\label{def:angle}
Given arbitrary element $b\in B_1^{\otimes L}$ where
number of \fbox{2} in $b$ is equal to or less than that of \fbox{1}.
Then the {\it angle variable} $[\iota (J)+d]\in\mathcal{J}(\mu)$
corresponding to $b$ is defined by the following procedure.
\begin{enumerate}
\item
Find integer $d$ and highest path $b_+$
such that $b=\bar{T}_1^d(b_+)$.
\item
Apply the bijection $\phi$ and obtain
$\phi (b_+)=\left((1^L),(\mu,J)\right)$.
\item
Extend the rigging to obtain $\iota (J)+d\in\bar{\mathcal{J}}(\mu)$.
\item
Take equivalent class with respect to $\simeq$ and
obtain $[\iota (J)+d]\in\mathcal{J}(\mu)$.
\hfill$\square$
\end{enumerate}
\end{definition}
\begin{remark}
The procedure in Definition \ref{def:angle}
uniquely determines the angle variable,
despite the non-uniqueness of $(d,b_+)$ in Step 1.
\hfill$\square$
\end{remark}

The following is the main theorem of \cite{KTT}.
\begin{theorem}\label{th:KTT}
Let the angle variable corresponding to $b$ be
$\left( J^{(k)}_i\right)_{i\in\mathbb{Z},k\in H}$.
Then the angle variable corresponding to $\bar{T}_l(b)$ is
$\left( J^{(k)}_i+\min (k,l)
\right)_{i\in\mathbb{Z},k\in H}$.
\hfill$\square$
\end{theorem}

\subsection{Ultradiscrete Riemann theta function}
In this section, we assume that all solitons have distinct lengths
$\mu_1<\mu_2<\cdots <\mu_g$.
More general case including solitons with same length can be
treated similarly (see \cite{KS2}).
Define ultradiscrete Riemann theta function
as follows:
\begin{equation}\label{eq:urt}
\begin{split}
\Theta({\bf z}) &= \lim_{\epsilon \rightarrow +0}
\epsilon\log\left( \sum_{{\bf n} \in \mathbb{Z}^g}
\exp\Bigl(-\frac{{}^t{\bf n}A{\bf n}/2+{}^t{\bf n}{\bf z}}{\epsilon}
\Bigr)\right)\\
&= -\min_{{\bf n} \in \mathbb{Z}^g}
\{{}^t{\bf n}A{\bf n}/2+{}^t{\bf n}{\bf z}\}.
\end{split}
\end{equation}
Here $A$
is the symmetric positive definite $g \times g$ integer matrix
appearing in the string center equation \cite{KN}:
\begin{equation}
(A)_{i,j} = \delta_{i,j}p_{\mu_i} + 2\min(\mu_i,\mu_j).
\end{equation}
Here $p_i$ is the vacancy number.

We introduce the vectors
\begin{eqnarray}
{\bf h}_l&=&(\min (\mu_i,l))_{i=1}^g\in\mathbb{Z}^g,\\
{\bf p}&=&(p_{\mu_i})_{i=1}^g\in\mathbb{Z}^g.
\end{eqnarray}
Again, consider the highest path $b_+$
obtained by $b=\bar{T}_1^d(b_+)$.
Let the rigging corresponding to $b_+$ be
$\mathbf{J}=(J_i)_{i=1}^g$.
Then we define $\mathbf{I}=(J_i+d)_{i=1}^g$.
\begin{definition}
For $1\leq k\leq L$ and $r=0,1$,
we define the ultradiscrete tau function as follows:
\begin{equation}\label{eq:tautheta}
\tau_r(k) = \Theta\left(
{\bf I}-\frac{\bf p}{2}-k{\bf h}_1+r{\bf h}_\infty\right).
\end{equation}
\hfill$\square$
\end{definition}
\begin{theorem}
Under the above settings,
the state $p$ is expressed as 
$p=(1-x(1),x(1))\otimes \cdots \otimes (1-x(L),x(L))$,
where
\begin{equation}
x(k)=\tau_0(k)-\tau_0(k-1)-\tau_1(k)+\tau_1(k-1).
\end{equation}
\hfill$\square$
\end{theorem}
Proof of this assertion uses
Proposition \ref{prop:highest*highest}
and the main result of \cite{KSY}.
See \cite{KS1}.
Note that this result itself is independent to
Theorem \ref{th:KTT}.

Combining this and Theorem \ref{th:KTT},
we solve the initial value problem of the pBBS.

\appendix
\section{Tensor product of highest paths}
In the appendix, we clarify special property of tensor product
of highest paths
which gives the basis for the inverse scattering formalism
of pBBS.
To begin with, we recall
famous characterization of highest paths.
\begin{lemma}
The highest elements
$b$ ($\tilde{e}_1 b=0$)
of the form $b_1\otimes\cdots\otimes b_L\in B_1^{\otimes L}$
are characterized by the following Yamanouchi condition:
\begin{equation}
\#\left\{1\leq i\leq k\,\Bigl|\, b_i=\fbox{1}\right\}\geq
\#\left\{1\leq i\leq k\,\Bigl|\, b_i=\fbox{2}\right\}
\end{equation}
for all $1\leq k\leq L$.
\hfill$\square$
\end{lemma}

In order to prove the following assertion,
we need to look at the original combinatorial
description of the map $\phi$
in addition to Theorem \ref{th:main}.
For description of the combinatorial algorithm of $\phi$,
see, e.g., Appendix A of \cite{KTT} or
Appendix C of \cite{KSY}.
Here we summarize basic definitions which will be used
in the proof.
Consider the rigged configuration
$\left((\lambda_i)_{i=1}^L ,(l_j,I_j)_{j=1}^N\right)$.
We call $\lambda$ quantum space and $(l,I)$ configuration.
Then the vacancy number $p_k$ for $k>0$ is defined by
\begin{equation}\label{def:app:vac}
p_k:=\sum_{i=1}^L\min (k,\lambda_i)
-2\sum_{j=1}^N\min (k,l_j).
\end{equation}
The row $l_j$ is called singular if the corresponding
rigging $I_j$ is equal to the vacancy number $p_{l_j}$
for the row $l_j$,
i.e., $p_{l_j}=I_j$.
Finally, we call quantity $p_{l_j}-I_j$ corigging.
It is known that $p_{l_j}\geq I_j$ for
all $(l_j,I_j)$.

In the following, we have to consider paths of the form
$q\otimes r$ where $q$ is arbitrary highest path
and $r$ is highest path of the form $r\in B_1^{\otimes M}$.
The basic points of the combinatorial procedure of $\phi$
in this setting are the following.
First of all, recall that the combinatorial procedure proceeds
recursively from the left of path to the right.
So we assume that we have done the procedure on $q$
and we exclusively consider the combinatorial
procedure on $r$.
To be more precise,
suppose that we have constructed the rigged configuration
corresponding to
$q\otimes r_{[k]}$ where $r_{[k]}$ is the first
$k$ components of $r$.
Then we are going to construct the rigged configuration
corresponding to $q\otimes r_{[k+1]}$
according to \fbox{1} or \fbox{2} of $(k+1)$th factor of $r$.
In both cases, we add one row of length one
to the quantum space of the rigged configuration
corresponding to $q\otimes r_{[k]}$.
For \fbox{2}, we add one box to the longest
singular row of configuration or, if there is no singular row,
we add one row of length one to the configuration.
The riggings for $q\otimes r_{[k+1]}$ are the same as
those for $q\otimes r_{[k]}$ except the rigging of row
of configuration that is different from $q\otimes r_{[k]}$.
We set the latter rigging equal to the vacancy number
(computed with the data of the rigged configuration
for $q\otimes r_{[k+1]}$)
for the corresponding row.

\begin{proposition}\label{prop:highest*highest}
Given two highest paths $q$ and $r$ as follows:
\begin{align}
q\in&\, B_{\lambda_1}\otimes B_{\lambda_2}\otimes\cdots\otimes
  B_{\lambda_L},\\
r\in&\, B_{\mu_1}\otimes B_{\mu_2}\otimes\cdots\otimes B_{\mu_M}.
\end{align}
Suppose that their rigged configurations are
$\phi (q)=\left(\lambda ,(l,I)\right)$
and $\phi (r)=\left( \mu,(m,K)\right)$.
Then the rigged configuration of the highest path
$q\otimes r$ is given by
$\phi (q\otimes r)=(\lambda\cup\mu ,(l\cup m,I\cup K'))$,
where $K'=\left( K^{'(j)}_i\right)$ is given by
\begin{equation}
K^{'(j)}_i=K^{(j)}_i+\tilde{p}_j,\qquad
\tilde{p}_j:=\sum_{k_1}\min (j,\lambda_{k_1})-2\sum_{k_2}\min (j,l_{k_2}),
\end{equation}
and $(l\cup m,I\cup K')$ means the union of $(l,I)$ and
$(m,K')$ as multi-sets of rows assigned with rigging.
$\tilde{p}_j$ is the vacancy number for $(\lambda ,(l,I))$.
\end{proposition}
{\bf Proof.}
Special version ($\lambda_i=1$ and $\mu_i=1$ for all $i$) of this claim
is proved in Lemma C.1 of \cite{KTT},
omitting some of details.
We include here an alternative proof intended
to clarify why the highest weight condition is
necessary for this result.

Consider the path $q\otimes r\otimes\fbox{1}^{\,\otimes\Lambda}$
with $\Lambda\gg |\mu |$.
Recall the property of the energy function $H(b\otimes u_l)=0$
for arbitrary $b\in B_k$.
This means that entries of the local energy distribution
under $\fbox{1}^{\,\otimes\Lambda}$ are all 0.
Therefore the rigged configuration obtained by
Theorem \ref{th:main} is the same as that for
the path $q\otimes r$ except
the extra $(1^\Lambda)$ of the quantum space.
So we can always think about paths of the form
$q\otimes r\otimes\fbox{1}^{\,\otimes\Lambda}$ by putting the tail
$\fbox{1}^{\,\otimes\Lambda}$ on the right of a given path.
Consider the isomorphism $q\otimes r\otimes\fbox{1}^{\,\otimes\Lambda}
\simeq q\otimes p^\ast\otimes\left(\bigotimes_i u_{\mu_i}\right)$,
where $p^\ast\in B_1^{\otimes\Lambda}$
and $p^\ast$ is highest.
{}From Lemma 8.5 of \cite{KSS}, these two isomorphic paths
correspond to the same rigged configuration.
Therefore we can assume $r\in B_1^{\otimes M}$
without loss of generality.

\renewcommand{\thefootnote}{\fnsymbol{footnote}}

Recall the Yamanouchi condition on $r$ and
consider the combinatorial procedure of $\phi$.
Assume that we have finished $\phi$ on
$q$ part and we are going to
apply $\phi$ to $r$.
Since we are assuming $r\in B_1^{\otimes M}$, all
letters 1 contained in $r$
correspond to length 1 row of the quantum space.
Fix a row of length $l$ that was constructed
from $q$, and consider the
change of corrigings induced by $r$.
During the procedure $\phi$,
if the chosen row does not obtain new box,
then the corresponding rigging does not change.
In such situation, we only need to keep track of change of the
vacancy numbers by using Eq.(\ref{def:app:vac}).
Then letters \fbox{1} of $r$ increase its corigging by 1,
on the other hand, letters \fbox{2} of $r$ decrease the coriggings
at most 1.
Therefore the Yamanouchi condition means that,
after creating rows corresponding to $q$,
those rows never become singular
during $r$.

Hence we can assume that rows corresponding to $r$
are independent to that of $q$, thus we have $(l\cup m)$
as a configuration.
In terms of the local energy distribution,
this means that no solitons cross the boundary
between $q$ and $r$.
As to the riggings, those corresponding to $m$
are larger than $K$ by $\tilde{p}_j$.
This follows from direct calculation using Eq.(\ref{s:eq:rigging}).
\vspace{5mm}
\hfill$\square$

\noindent
{\bf Acknowledgements:}
The author is grateful to Prof. Masato Okado for
organizing a nice conference
``Expansion of Combinatorial Representation Theory"
(RIMS, Kyoto University, October 2007)
where the lecture was given
and to Prof. Atsuo Kuniba for collaboration.
He is a research fellow of the 
Japan Society for the Promotion of Science.

\bigskip

\begin{flushleft}
Reiho Sakamoto\bigskip\\
Department of Physics, Graduate School of Science,\\
The University of Tokyo\\
Hongo, Bunkyo-ku, Tokyo, 113-0033, Japan\\
e-mail: \texttt{reiho@spin.phys.s.u-tokyo.ac.jp}
\end{flushleft}
\end{document}